\documentclass[oneside, 12pt]{amsart}

\usepackage[margin=1in]{geometry}
\usepackage{amsthm}
\usepackage{enumitem}  
\usepackage[OT2,T1]{fontenc}
\usepackage{amscd, amssymb, enumitem, amsmath, mathrsfs, amsfonts, amsaddr}
\usepackage{url}
\usepackage{tikz}
\usepackage{fullpage}
\usepackage{microtype}
\usetikzlibrary{decorations.pathreplacing}
\usepackage[backref=page]{hyperref}
\usepackage{hyperref}
\usepackage[utf8]{inputenc}

\DeclareSymbolFont{cyrletters}{OT2}{wncyr}{m}{n}
\DeclareMathSymbol{\Sha}{\mathalpha}{cyrletters}{"58}

\newtheorem{theorem}{Theorem}[section]

\newtheorem{lemma}[theorem]{Lemma}
\newtheorem{proposition}[theorem]{Proposition}
\newtheorem*{proposition*}{Proposition}
\newtheorem{corollary}[theorem]{Corollary}

\theoremstyle{definition}

\newtheorem{conjecture}[theorem]{Conjecture}
\newtheorem{claim}[theorem]{Claim}
\newtheorem{remark}[theorem]{Remark}
\newtheorem{emptyremark}[theorem]{}
\newtheorem*{acknowledgement}{Acknowledgements}

\theoremstyle{remark}

\title{On a Conjecture of Agashe}
\author{Mentzelos Melistas}
\address{Department of Mathematics, University of Georgia, Athens, GA 30602,\\ ~~~~~~~~~~~~~~~~~~~~~~~~~~~~~~~~~~ email: mentzmel@gmail.com}
\date{\today. \small{2010 {\it Mathematics Subject Classification.} Primary 11G05;} Secondary 11G07, 11G40}

\begin{document}

\maketitle

\begin{abstract}
    Let $E/\mathbb{Q}$ be an optimal elliptic curve, $-D$ be a negative fundamental discriminant coprime to the conductor $N$ of $E/\mathbb{Q}$ and let $E^{-D}/\mathbb{Q}$ be the twist of $E/\mathbb{Q}$ by $-D$. A conjecture of Agashe predicts that if $E^{-D}/\mathbb{Q}$ has analytic rank $0$, then the square of the order of the torsion subgroup of $E^{-D}/\mathbb{Q}$ divides the product of the order of the Shafarevich-Tate group of $E^{-D}/\mathbb{Q}$ and the orders of the arithmetic component groups of $E^{-D}/\mathbb{Q}$, up to a power of $2$. This conjecture can be viewed as evidence for the second part of the Birch and Swinnerton-Dyer conjecture for elliptic curves of analytic rank zero. We provide a proof of a slightly more general statement without using the optimality hypothesis.

\end{abstract}

\section{Introduction}
Let $E/\mathbb{Q}$ be an optimal elliptic curve and let $-D$ be a negative fundamental discriminant, i.e., $-D$ is negative and is the discriminant of a quadratic number field. Assume that $D$ is coprime to the conductor $N$ of $E/\mathbb{Q}$ and let $E^{-D}/\mathbb{Q}$ denote the twist of $E/\mathbb{Q}$ by $-D$. Let $L(E^{-D},s)$ be the $L$-function of $E^{-D}/\mathbb{Q}$ and assume that $L(E^{-D},1)\neq 0$, i.e., that $E^{-D}/\mathbb{Q}$
has analytic rank $0$. In \cite{aga}, under some additional mild hypotheses, Agashe proved that $\frac{L(E^{-D},1)}{\Omega(E^{-D})} \in \mathbb{Z}[\frac{1}{2}]$. The second part of the Birch and Swinnerton-Dyer conjecture for elliptic curves of rank $0$ predicts that $$\frac{L(E^{-D},1)}{\Omega(E^{-D})}= \frac{|\Sha(E^{-D}/\mathbb{Q})| \cdot \prod_{p}c_p(E^{-D})}{|E^{-D}(\mathbb{Q})|^2},$$ see Conjecture $1.2$ of \cite{aga}. Here $\Sha(E^{-D}/\mathbb{Q})$ denotes the Shafarevich-Tate group of $E^{-D}/\mathbb{Q}$ and $c_p(E^{-D})$ is the order of the arithmetic component group of $E^{-D}/\mathbb{Q}$ at $p$, also called the Tamagawa number of $E^{-D}/\mathbb{Q}$ at $p$. This, along with results of \cite{pal}, led Agashe, in Section $5.3$ of the appendix of \cite{pal}, to propose the following conjecture.

\begin{conjecture}\label{agasheconjecture}
Let $E/\mathbb{Q}$ be an optimal elliptic curve of conductor $N$ and let $-D$ be a negative fundamental discriminant such that $D$ is coprime to $N$. Let $E^{-D}/\mathbb{Q}$ denote the twist of $E/\mathbb{Q}$ by $-D$. Suppose that $L(E^{-D},1) \neq 0$. Then
\vspace{-.1cm}
\begin{align*}
 |E^{-D}(\mathbb{Q})|^2 \quad \text{divides} \quad |\Sha(E^{-D}/\mathbb{Q})| \cdot \prod_{p|N}c_p(E^{-D}), \quad \text{up to a power of 2.} 
\end{align*}
\end{conjecture}

Theorem \ref{maintheorem} below, which is the main theorem of this paper, implies Conjecture \ref{agasheconjecture} (see Corollary \ref{corollaryagasheconjecture} for the proof of this implication). We note that Conjecture \ref{agasheconjecture} and, hence, our theorem can be seen as further evidence for the second part of the Birch and Swinnerton-Dyer conjecture in the analytic rank zero case.

\begin{theorem}\label{maintheorem} {\bf (see Corollary \ref{corollarytorsion} and Theorem \ref{prop:dpowerof3}) }
Let $E/\mathbb{Q}$ be an elliptic curve of conductor $N$ and let $d \neq \pm 1$ be a square-free integer such that $d$ is coprime to $N$.
\begin{enumerate}
     \item If $d\neq \pm3$, then the torsion subgroup of $E^d(\mathbb{Q})$ contains only points of order dividing $8$.
     \end{enumerate}
\noindent  Assume also that $L(E^d,1) \neq 0$. Then the groups $E^d(\mathbb{Q})$ and $|\Sha(E^d/\mathbb{Q})|$ are finite and, up
to a power of $2$, the following conditions hold.
\begin{enumerate}[label=(\roman*),resume]
     \item If $d \neq 3$, then $
 |E^d(\mathbb{Q})|^2 \quad \text{divides} \quad |\Sha(E^d/\mathbb{Q})| \cdot \displaystyle\prod_{p|N}c_p(E^d).$ 
     \item If $d=3$, then $ 
 |E^d(\mathbb{Q})|^2 \quad \text{divides} \quad |\Sha(E^d/\mathbb{Q})| \cdot \displaystyle\prod_{p|2N}c_p(E^d).$
\end{enumerate}
\end{theorem}

\begin{remark}
If moreover $d \neq -3$, then Part $(i)$ of Theorem \ref{maintheorem} implies Part $(ii)$. We also note that the case $d = \pm 3$ is by far
the most difficult case to treat in the proof.
\end{remark}

\begin{remark}
Clearly, the analogue of Part (i) of Theorem \ref{maintheorem} with $d=1$ does not hold, since there exist elliptic curves $E/\mathbb{Q}$ of rank $0$ which have $\mathbb{Q}$-rational torsion points of order 3, 5, or 7.
Consider now the elliptic curve with LMFDB \cite{lmfdb} label 176.a2. Then $E^{-1}/\mathbb{Q}$, the twist of $E/\mathbb{Q}$ by $-1$, is the elliptic curve with LMFDB label 44.a2 and $E^{-1}(\mathbb{Q})=\mathbb{Z}/3\mathbb{Z}$. This proves that Part $(i)$ of Theorem \ref{maintheorem} is false if $d=-1$. 

Moreover, $|\Sha(E^{-1}/\mathbb{Q})|=1$ and $\prod_{p} c_p(E^{-1})=c_2(E^{-1})=3$. Therefore, $|E^{-1}(\mathbb{Q})|^2=9$ while $|\Sha(E^{-1}/\mathbb{Q})|\prod_{p} c_p(E^{-1})=3$. This proves that the analogues of Theorem \ref{maintheorem}, Parts $(ii), (iii)$ are false when $d= -1$.
\end{remark}

The proof of Part $(i)$ of Theorem \ref{maintheorem} is postponed to Corollary \ref{corollarytorsion} and the proof of Parts $(ii)$ and $(iii)$ is completed in Theorem \ref{prop:dpowerof3}. To prove Parts $(ii)$ and $(iii)$ we prove Theorem \ref{3torsionshaandtamagawa}, which has applications to another consequence of the Birch and Swinnerton-Dyer conjecture, discovered by Agashe and Stein in \cite{agashestein} (see Remark \ref{remarkagashestein}).

This paper is organized as follows. First, we prove that if $d \neq \pm 1, \pm 3$, then $E^d(\mathbb{Q})$ only contains points of order a power of $2$. This explains some observations made by Agashe, in the case where $d$ is a negative fundamental discriminant, which can be found on page $8$ of \cite{aga}. Next, we prove that for $d=-3$ or $3$, $E^d(\mathbb{Q})$ cannot contain points of order $5$ or $7$. All this is done in section \ref{section2}. Finally, in section \ref{section3} we consider in detail the case where $d=\pm 3$ and $E^d(\mathbb{Q})$ contains a point of order $3$.
 
 \begin{acknowledgement}
This work is part of the author's doctoral dissertation at the University of Georgia. The author would like to thank his advisor Dino Lorenzini for valuable help during the preparation of this work as well as for many useful suggestions on improving the exposition of this manuscript. The author thanks the anonymous referee for many insightful comments.
 \end{acknowledgement}
 
 \section{Restrictions on the torsion subgroup}\label{section2}

Let $E/\mathbb{Q}$ be an elliptic curve and let $d$ be any non-zero integer. Recall that the twist of $E/\mathbb{Q}$ by $d$, which is denoted by $E^d/\mathbb{Q}$, is an elliptic curve over $\mathbb{Q}$ that becomes isomorphic over $\mathbb{Q}(\sqrt{d})$ to $E_{\mathbb{Q}(\sqrt{d})}/\mathbb{Q}(\sqrt{d})$ but which is not in general isomorphic to $E/\mathbb{Q}$ over $\mathbb{Q}$. More explicitly, if $E/\mathbb{Q}$ is given by a Weierstrass equation of the form $y^2=x^3+ax+b$, then $E^d/\mathbb{Q}$ is given by $dy^2=x^3+ax+b$ or, after making a change of variables, by $y^2=x^3+d^2ax+d^3b$.

In this section, we obtain results that put restrictions on the rational torsion subgroup of $E^d/\mathbb{Q}$ provided that $d$ and the conductor of $E/\mathbb{Q}$ are coprime.

\begin{proposition}\label{no5or7torsion}
Let $E/\mathbb{Q}$ be an elliptic curve defined over $\mathbb{Q}$ of conductor $N$, and let $d \neq \pm 1$ be any square-free integer that is coprime to $N$. Then $E^d(\mathbb{Q})$ cannot contain a rational point of order $5$ or $7$.
\end{proposition}
\begin{proof}
We prove the proposition by contradiction. Let $\ell$ be equal to $5$ or $7$, assume that $E^d(\mathbb{Q})$ contains a point of order $\ell$, and that $d \neq \pm 1$. Since $E^d/\mathbb{Q}$ has a $\mathbb{Q}$-rational point of order $\ell$, Proposition $8.4$ in Chapter XV of \cite{mfl} implies that $E^d/\mathbb{Q}$ is semi-stable away from $\ell$.

Assume that there exists an odd prime $p$ with $p \mid d$. Then $p \nmid N$ because $d$ and $N$ are coprime. Therefore, Proposition $1$ of \cite{com} implies that $E^d/\mathbb{Q}$ has reduction type I$_0^*$ modulo $p$.  As a result, since $E^d/\mathbb{Q}$ is semi-stable away from $\ell$, we obtain that $p$ can only be equal to $\ell$. However, Proposition \ref{mentzelospaper5,7} below implies that $E^d/\mathbb{Q}$ cannot have reduction of type I$_0^*$ modulo $\ell$. Therefore, $d= \pm 2$ since $d$ is square-free. 

Assume now that $d=-2$ or $2$. We will arrive at a contradiction. Since we assume that $d$ and $N$ are coprime, we get that $2 \nmid N$. Since $E/\mathbb{Q}$ has good reduction modulo $2$, by Table II of \cite{com} (see also the hypotheses for this table at the bottom of page 58 of \cite{com}) we obtain that $E^d/\mathbb{Q}$ has reduction of type I$_8^*$ or II modulo $2$. In \cite{com}, good reduction is denoted as usual with the symbol I$_0$, and our twist $E^d$ is denoted by $E^{\chi}$. Note that Table II is independent of the existence of a torsion point. We now use our assumption that $E^d$ has a $\mathbb{Q}$-rational point of order $\ell$. Since $E^d$ has a $\mathbb{Q}$-rational point of order $\ell$, it has semi-stable reduction modulo $2$, which is a contradiction. This proves that $E^d(\mathbb{Q})$ cannot contain a point of order $\ell$.
\end{proof}

A proof of the following proposition can be found in the forthcoming paper \cite{mentzeloskodairandtorsion}.
\begin{proposition}\label{mentzelospaper5,7}
Let $E/\mathbb{Q}$ be an elliptic curve  with a $\mathbb{Q}$-rational point of order $p$, where $p$ is a prime number. Assume that $E/\mathbb{Q}$ has additive reduction modulo $p$. Then $p \leq 7$, and 
 \begin{enumerate}[label=(\roman*),topsep=2pt]\label{5and7torsionreduction}
    \itemsep0em 
    \item If $p=5$, then $E/\mathbb{Q}$ can only have reduction of type II or III modulo $5$.
    \item If $p=7$, then $E/\mathbb{Q}$ can only have reduction type II modulo $7$.
 \end{enumerate}
\end{proposition}

We now investigate whether $E^d(\mathbb{Q})$ can contain a point of order $3$. Let $E/\mathbb{Q}$ be an elliptic curve with a $\mathbb{Q}$-rational point of order $3$. Then, by translating that point to $(0,0)$ and performing a change of variables if necessary (see Remark $2.2$ in Section $4.2$ of \cite{hus}) we obtain a Weierstrass equation for $E/\mathbb{Q}$ of the form $$y^2+cxy+dy=x^3,$$ with $c,d \in \mathbb{Q}$. If $u \in \mathbb{Z}$, then the transformation $(x,y) \rightarrow{(\frac{x}{u^2}, \frac{y}{u^3})}$ gives a new Weierstrass equation of the same form with $c$ replaced by $uc$ and with $d$ replaced by $u^3d$ (see page 185 of \cite{aec}). Therefore, by picking $u$ to be a large power of the product of the denominators of $c,d$ (if any), we can arrange that $c,d \in \mathbb{Z}$. Moreover, by applying the transformation $(x,y) \rightarrow{(x, -y)}$ if necessary, we can arrange that $d > 0$. We now show that we can find a new Weierstrass equation of the above form with coefficients, which we will call $a$ and $b$ below, that also have the property that for every prime $q$ either $q \nmid a$ or $q^3 \nmid b$. If there is no prime $q$ such that $q \mid c$ and $q^3 \mid d$, then set $a=c$ and $b=d$. Otherwise, let $q_1,...,q_s$ be the set of primes such that $q_i \mid c$ and ${q_i}^3 \mid d$, and let $n_i= \mathrm{min}\{ \mathrm{ord}_{q_i}(c), \lfloor \frac{\mathrm{ord}_{q_i}(d)}{3} \rfloor \} $, where $\lfloor -  \rfloor$ is the floor function. Let $u=\prod_{i=1}^{s} q_i^{n_i}$. By applying the transformation $(x,y) \rightarrow{(u^2x, u^3y)}$, we obtain a new Weierstrass equation of the form $$y^2+\frac{c}{u}xy+\frac{d}{u^3}y=x^3.$$ Setting $a=\frac{c}{u}$ and $b=\frac{d}{u^3}$ we see that $a,b \in \mathbb{Z}$, $b >0$, and for every prime  $q$ either $q \nmid a$ or $q^3 \nmid b$. Therefore, we have proved that we can choose a Weierstrass equation for $E/\mathbb{Q}$ of the form
\begin{align*}\label{eq:3torsion}
    y^2+axy+by=x^3, \tag{$2.3$}
\end{align*}
where $a,b$ are integers, $b>0$, and for every prime $q$ either $q \nmid a$ or $q^3 \nmid b$. Also, we must have $a^3-27b \neq 0$, since the discriminant of Equation $($\ref{eq:3torsion}$)$ is $$\Delta=b^3(a^3-27b) \quad \text{and we also have} \quad c_4=a(a^3-24b).$$
\stepcounter{theorem}
\begin{proposition}\label{prop:3torsionreduction} (See Proposition $3.5$ and Lemma $3.6$ of \cite{koz})
Let $E/\mathbb{Q}$ be an elliptic curve given by Equation $($\ref{eq:3torsion}$)$. Write $D:=a^3-27b$ and let $p$ be any prime (note that either $p \nmid a$ or $p^3 \nmid b$). Then the reduction of $E/\mathbb{Q}$ modulo $p$ is determined as follows: 
\begin{enumerate}
    \item $\text{If} \; 3ord_p(a) \leq ord_p(b)$ :  $\begin{cases} 
      3ord_p(a)< ord_p(b) \rightarrow \text{split I}_{3ord_p(b)}, \quad c_p(E)=3ord_p(b) \\
      3ord_p(a)= ord_p(b) \rightarrow \begin{cases} 
      ord_p(D)>0 \rightarrow \text{I}_{ord_p(D)}\\
      ord_p(D)=0 \rightarrow \text{Good reduction I}_0
   \end{cases}
\end{cases}$

    \item $\text{If} \; 3ord_p(a) > ord_p(b)$ :  $\begin{cases} 
      ord_p(b) = 0 \rightarrow \begin{cases} 
      p=3 \rightarrow \text{Go to} \; (iii) \\
      p \neq 3 \rightarrow \text{Good reduction I}_0
   \end{cases} \\
      ord_p(b) = 1  \rightarrow \text{IV, } c_p(E)=3\\
       ord_p(b) = 2 \rightarrow \text{IV}^*,\text{ } c_p(E)=3.
\end{cases}$

   \item If $p=3$ and $ord_p(a)>0=ord_p(b)$: $ord_p(D)=$ $\begin{cases}
3 \rightarrow \text{II or III}\\
4 \rightarrow \text{II}\\
5 \rightarrow \text{IV}\\
n \rightarrow \text{I}_{n-6}^* \text{  , for  } n \geq 6.
\end{cases}$

\end{enumerate}
\end{proposition}

\begin{proposition}\label{no3torsiondnotequaltopm3}
Let $E/\mathbb{Q}$ be an elliptic curve over $\mathbb{Q}$ of conductor $N$ and let $d \neq \pm1, \pm 3$ be a square-free integer coprime to $N$. Then $E^d(\mathbb{Q})$ cannot contain a point of order $3$.
\end{proposition}
\begin{proof}
The proof is by contradiction. Assume that $E^d(\mathbb{Q})$ contains a point of order $3$ and that $d$ has a prime divisor $p \neq 3$. Assume first that $p \neq 2$. Since $p \nmid N$, Proposition $1$ of \cite{com} implies that $E^d/\mathbb{Q}$ has reduction of type I$_0^*$ modulo $p$. However, this is impossible because by Proposition \ref{prop:3torsionreduction}, $E^d(\mathbb{Q})$ cannot have reduction of type I$_0^*$ modulo $p$ for $p \neq 3$. Assume now that $p=2$. Since $2 \nmid N$, we obtain that $E/\mathbb{Q}$ has good reduction modulo $2$ and, hence, by Table II  of \cite{com} (see also the hypotheses for this table at the bottom of page 58 of \cite{com}), we get that $E^d/\mathbb{Q}$ has reduction of type I$_8^*$ or II modulo $2$. However, $E^d/\mathbb{Q}$ cannot have reduction of type I$_8^*$ or II modulo $2$ by Proposition \ref{prop:3torsionreduction}. This implies that $2 \nmid d$. Therefore, $d$ can only be divisible by $3$ so $d= \pm 1$ or $\pm3$.
\end{proof}

\begin{corollary}\label{corollarytorsion}
Let $E/\mathbb{Q}$ be an elliptic curve over $\mathbb{Q}$ of conductor $N$ and let $d \neq \pm1$ be a square-free integer such that $d$ is coprime to $N$. Then
 \begin{enumerate}[label=(\roman*),topsep=2pt]
    \itemsep0em 
    \item $E^d(\mathbb{Q})_{tors}$ contains only points of order $2^{\alpha}3^{\beta}$ where $\alpha, \beta \geq 0$.
    \item If $d \neq \pm3$, then $E^d(\mathbb{Q})_{tors}$ contains only points of order dividing $8$.
 \end{enumerate}
\end{corollary}
\begin{proof}
By a Theorem of Mazur (see \cite{maz} Theorem $(8)$) the only primes that can divide $|E^d(\mathbb{Q})_{\mathrm{tors}}|$ are $2,3,5$ and $7$. Therefore, Proposition \ref{no5or7torsion} implies Part $(i)$. Moreover, Propositions \ref{no5or7torsion} and \ref{no3torsiondnotequaltopm3} imply that if $d \neq \pm 1, \pm 3$, then $E^d(\mathbb{Q})_{tors}$ contains only points of order dividing $2^{\alpha}$, where $\alpha \geq 0$. Finally, Part $(ii)$ follows by applying Mazur's Theorem (see \cite{maz} Theorem $(8)$).
\end{proof}

\section{Elliptic curves over $\mathbb{Q}$ with a rational point of order 3}\label{section3}

In this section, we prove Parts $(ii)$ and $(iii)$ of Theorem \ref{maintheorem} in Theorem \ref{prop:dpowerof3}. To achieve this goal, we will use Theorem \ref{3torsionshaandtamagawa} below, which might be of independent interest (see Remark \ref{remarkagashestein}).

\begin{theorem}\label{3torsionshaandtamagawa}
Let $E/\mathbb{Q}$ be an elliptic curve with a $\mathbb{Q}$-rational point of order $3$. Assume that the analytic rank of $E/\mathbb{Q}$ is $0$ and that $E/\mathbb{Q}$ has reduction of type I$_n^*$ modulo $3$, for some $n \geq 0$.
\begin{enumerate}[label=(\alph*)]
    \item If $E/\mathbb{Q}$ has semi-stable reduction away from $3$, then $9 \mid |\Sha(E/\mathbb{Q})| \cdot \prod_{p}c_p(E)$.
    \item If $E/\mathbb{Q}$ has more than two places of additive reduction, then $9 \mid \prod_{p}c_p(E)$.
    \item If $E/\mathbb{Q}$ has exactly two places of additive reduction, then $9 \mid |\Sha(E/\mathbb{Q})| \cdot \prod_{p}c_p(E)$.
\end{enumerate}
\end{theorem}

\begin{proof}[\it Proof of Theorem \ref{3torsionshaandtamagawa}]
Let $j_E$ be the $j$-invariant of $E/\mathbb{Q}$. Since some of our arguments below do not work for $j_E=0$ or $1728$, we first handle these cases with the following claim.
\begin{claim}\label{claimjinvariants}
 Let $E/\mathbb{Q}$ be an elliptic curve with $j$-invariant $0$ or $1728$ that has a $\mathbb{Q}$-rational point of order $3$. Then $E/\mathbb{Q}$ cannot have reduction of type I$_n^*$ modulo $3$, for any $n \geq 0$.
\end{claim}
\begin{proof}
Let $E/\mathbb{Q}$ be an elliptic curve with $j$-invariant $j_E=0$. Since, $j_E=\frac{{c_4}^3}{\Delta}$ we obtain that $c_4=0$. Tableau II of \cite{pap} implies that $E/\mathbb{Q}$ can only have reduction type II, II$^*$, III, III$^*$, IV, or IV$^*$ modulo $3$. Therefore, $E/\mathbb{Q}$ cannot have reduction type I$_n^*$ modulo $3$, for any $n \geq 0$, as needed.

Let $E/\mathbb{Q}$ be an elliptic curve with $j$-invariant $j_E=1728$ and with a $\mathbb{Q}$-rational point of order $3$. Since $E/\mathbb{Q}$ has a $\mathbb{Q}$-rational point of order $3$, it has a Weierstrass equation of the form $($\ref{eq:3torsion}$)$. Assume that $E/\mathbb{Q}$ has reduction type I$_n^*$ modulo $3$, for some $n \geq 0$, and we will find a contradiction. Proposition \ref{prop:3torsionreduction}, Parts $(i)$ and $(ii)$ imply that $3 \nmid b$ and $3 \mid a$. Since $j_E=\frac{c_4^3}{\Delta}$ and $j_E=1728$, we obtain that $c_4^3=1728\Delta$. Writing $c_4$ and $\Delta$ in terms of $a$ and $b$, we get that $(a(a^3-24b))^3=1728b^3(a^3-27b)$ and, hence, $$3{\rm ord}_3(a)+3{\rm ord}_3(a^3-24b)={\rm ord}_3(1728)+3{\rm ord}_3(b)+{\rm ord}_3(a^3-27b).$$ Since $3 \mid a$ and $3 \nmid b$, we get that ord$_3(a^3-24b)=1$. Therefore, $$3{\rm ord}_3(a)+3=3+0+{\rm ord}_3(a^3-27b),$$ which yields
$$3{\rm ord}_3(a)={\rm ord}_3(a^3-27b).$$
If ${\rm ord}_3(a)=1$, then ${\rm ord}_3(a^3-27b)=3$, which is a contradiction because Proposition \ref{prop:3torsionreduction}, Part $(iii)$ implies that $E/\mathbb{Q}$ has reduction type II or III. If ${\rm ord}_3(a) \geq 2$, then ${\rm ord}_3(a^3-27b)=3$ because ${\rm ord}_3(b)=0$, and this is again a contradiction. This proves our claim.
\end{proof}
Claim \ref{claimjinvariants} implies that to prove our theorem, it is enough to assume that $j_E \neq 0, 1728$. We assume from now on that $j_E \neq 0, 1728$.

Since $E/\mathbb{Q}$ has a $\mathbb{Q}$-rational point of order $3$, it can be given by an equation of the form $($\ref{eq:3torsion}$)$ where $(0,0)$ is a $\mathbb{Q}$-rational point of order $3$. Let $\widehat{E}:=E/<(0,0)>$ and let $\phi:  E \rightarrow{} \widehat{E}$ be the associated $3$-isogeny. We denote by $\hat{\phi}$ the dual isogeny.

We recall now some important facts in preparation for the proof of Lemma \ref{tamagawagreaterorequalto2} below. The following formula is due to Cassels (see Theorem $1$ of \cite{ks} or Theorem $2.1$ of \cite{koz})
\begin{align}\label{eq:selmerandtamagawa}
 \frac{|\mathrm{Sel}^{(\phi)}(E/\mathbb{Q})|}{|\mathrm{Sel}^{(\hat{\phi})}(\widehat{E}/\mathbb{Q})|}=\frac{|E[\phi](\mathbb{Q})|\Omega(\widehat{E})\prod_{p}c_p(\widehat{E})}{|\widehat{E}[\hat{\phi}](\mathbb{Q})|\Omega(E)\prod_{p}c_p(E)}.  \tag{3.3}
\end{align}
Here $\mathrm{Sel}^{(\phi)}(E/\mathbb{Q})$ and  $\mathrm{Sel}^{(\hat{\phi})}(\widehat{E}/\mathbb{Q})$ are the $\phi$- and $\hat{\phi}$-Selmer groups, respectively. Moreover, $\Omega(E):=\displaystyle\int_{E(\mathbb{R})}|\omega_{min}|$ and $\Omega(\widehat{E}):=\displaystyle\int_{\widehat{E}(\mathbb{R})}|\widehat{\omega}_{min}|$, where $\omega_{min}$ and $\widehat{\omega}_{min}$ are two minimal invariant differentials on $E/\mathbb{Q}$ and $\widehat{E}/\mathbb{Q}$ respectively (see Section III.1 of \cite{aec} and page 451 of \cite{aec}). The interested reader can consult section X.$4$ of \cite{aec} concerning the definitions of Selmer groups of elliptic curves as well as relevant background.

Since $E[\phi](\mathbb{Q}) \cong \mathbb{Z}/3\mathbb{Z}$, we get that $\widehat{E}[\hat{\phi}](\mathbb{Q})$ is trivial by the Weil pairing (see \cite{aec} exercise III.3.15 on page $108$). Hence,
\stepcounter{theorem}
\begin{align}\label{kernelofisogeny}
 \frac{|E[\phi](\mathbb{Q})|}{|\widehat{E}[\hat{\phi}](\mathbb{Q})|}=3. \tag{3.4}
\end{align}

We know that $\phi^*\widehat{\omega}_{min}=\lambda \omega_{min}$, for some $\lambda \in \mathbb{Z}$ (see the second paragraph of page $284$ of \cite{vatsal}). Moreover, the fact that $\hat{\phi}\circ \phi=[3]$ implies that $|\lambda|=1$ or $3$.
Therefore, $|\frac{\omega_{min}}{\phi^*\widehat{\omega}_{min}}|$ is equal to $1$ or $3^{-1}$. By Lemma $7.4$ of \cite{dd},
\begin{align*}
     \frac{\Omega(E)}{\Omega(\widehat{E})}=\frac{|{\rm ker} (\phi :E(\mathbb{R}) \longrightarrow \widehat{E}(\mathbb{R}))|}{|{\rm coker} (\phi :E(\mathbb{R}) \longrightarrow \widehat{E}(\mathbb{R}) )|}\cdot|\frac{\omega_{min}}{\phi^*\widehat{\omega}_{min}}|
\end{align*}
Since ${\rm ker} (\phi) \subset E(\mathbb{R})$ and $\phi$ has degree $3$, Proposition $7.6$ of \cite{dd} implies that $\frac{|{\rm ker} \phi :E(\mathbb{R}) \longrightarrow \widehat{E}(\mathbb{R})|}{|{\rm coker} (\phi :E(\mathbb{R}) \longrightarrow \widehat{E}(\mathbb{R}) )|}=3$. Therefore, since $|\frac{\omega_{min}}{\phi^*\widehat{\omega}_{min}}|$ is equal to $1$ or $3^{-1}$,
\stepcounter{theorem}
\begin{align}\label{Omega}
 \frac{\Omega(E)}{\Omega(\widehat{E})}=1 \; \text{or} \; 3. \tag{3.5}
\end{align}
\stepcounter{theorem}
The following lemma will be used repeatedly in what follows. In pursuing this idea, we were inspired by work of Byeon, Kim, and Yhee in \cite{bky}.
\begin{lemma}\label{tamagawagreaterorequalto2}
Let $E/\mathbb{Q}$ be an elliptic curve given by a Weierstrass equation of the form $($\ref{eq:3torsion}$)$ and assume that $E/\mathbb{Q}$ has analytic rank $0$. Let $\widehat{E}:=E/<(0,0)>$ and let $\phi:  E \rightarrow{} \widehat{E}$ be the associated $3$-isogeny. 
\begin{enumerate}
    \item If $\mathrm{ord}_3(\frac{\prod_{p}c_p(\widehat{E})}{\prod_{p}c_p(E)})\geq 2$ or $\rm{dim}_{\mathbb{F}_3}(\mathrm{Sel}^{(\phi)}(E/\mathbb{Q})) \geq 2$, then $9$ divides $|\Sha(E/\mathbb{Q})|$.
    \item If $\frac{|\mathrm{Sel}^{(\phi)}(E/\mathbb{Q})|}{|\mathrm{Sel}^{(\hat{\phi})}(\widehat{E}/\mathbb{Q})|}=1$, then $9$ divides $|\Sha(E/\mathbb{Q})|$.
\end{enumerate}

\end{lemma}
\begin{proof}\text{}
Since the analytic rank of $E/\mathbb{Q}$ is zero, work of Gross and Zagier, on heights of Heegner points \cite{grosszagierpaper}, as well as work of Kolyvagin, on Euler systems \cite{kolyvagineulersystems}, imply that  $E/\mathbb{Q}$ has (algebraic) rank $0$ and that $\Sha(E/\mathbb{Q})$ is finite (see Theorem $3.22$ of \cite{darmonmodularellipticcurves} for a sketch of the proof).
Thus the statements of the lemma make sense.

{\it Proof of $(i)$:} First we prove that if dim$_{\mathbb{F}_3}(\mathrm{Sel}^{(\phi)}(E/\mathbb{Q})) \geq 2$, then $9$ divides $|\Sha(E/\mathbb{Q})|$. There is a short exact sequence $$0 \longrightarrow{} E(\mathbb{Q})/3E(\mathbb{Q}) \longrightarrow{} \mathrm{Sel}^{3}(E/\mathbb{Q}) \longrightarrow{} \Sha(E/\mathbb{Q})[3] \longrightarrow{} 0$$ (see Theorem X.$4.2$ of \cite{aec}).
Therefore, using the fact that $E(\mathbb{Q})/3E(\mathbb{Q}) \cong \mathbb{Z}/3\mathbb{Z}$, since $E/\mathbb{Q}$ has rank $0$ and $E[3](\mathbb{Q}) \cong \mathbb{Z}/3\mathbb{Z}$, we see that dim$_{\mathbb{F}_3}(\mathrm{Sel}^{3}(E/\mathbb{Q})) \geq 2$ implies that $\Sha(E/\mathbb{Q})[3]$ has a nontrivial element. Since the order of $\Sha(E/\mathbb{Q})$ is a square (see \cite{aec} Corollary $17.2.1$) if $\Sha(E/\mathbb{Q})[3]$ has a nontrivial element, then $9$ divides $|\Sha(E/\mathbb{Q})|$. By Corollary 1, section 2 of \cite{ks} and the fact that $\widehat{E}[\hat{\phi}](\mathbb{Q})$ is trivial, we obtain that the natural map $\mathrm{Sel}^{(\phi)}(E/\mathbb{Q}) \longrightarrow \mathrm{Sel}^{3}(E/\mathbb{Q})$ is injective. Therefore, if dim$_{\mathbb{F}_3}(\mathrm{Sel}^{(\phi)}(E/\mathbb{Q})) \geq 2$, then $9$ divides $|\Sha(E/\mathbb{Q})|$. 

Finally, if  $\text{ord}_3(\frac{\prod_{p}c_p(\widehat{E})}{\prod_{p}c_p(E)})\geq 2$, then Equation $($\ref{eq:selmerandtamagawa}$)$, combined with $($\ref{kernelofisogeny}$)$ and $($\ref{Omega}$)$, implies that dim$_{\mathbb{F}_3}(\mathrm{Sel}^{(\phi)}(E/\mathbb{Q})) \geq 2$ so by the first part of the proof, we find that $9$ divides $ |\Sha(E/\mathbb{Q})|$. This concludes the proof of Part $(i)$.

{\it Proof of $(ii)$:}
By Corollary 1 of section $2$ of \cite{ks} and the fact that $\widehat{E}[\hat{\phi}](\mathbb{Q})$ is trivial there is an exact sequence:
\stepcounter{theorem}
\begin{align*}\label{exactsequencephiselmergroup}
    0 \rightarrow \mathrm{Sel}^{(\phi)}(E/\mathbb{Q}) \rightarrow \mathrm{Sel}^3 (E/\mathbb{Q}) \rightarrow \mathrm{Sel}^{(\hat{\phi})}(\widehat{E}/\mathbb{Q}) \rightarrow \Sha(\widehat{E}/\mathbb{Q})[\hat{\phi}]/\phi(\Sha(E/\mathbb{Q})[3]) \rightarrow 0. \tag{3.7}
\end{align*}
Moreover, by Theorem $4.3$ of \cite{shnidmanquadraticrm}, the Cassels-Tate pairing on $\Sha(\widehat{E}/\mathbb{Q})$ restricts to a non-degenerate alternating pairing on $\Sha(\widehat{E}/\mathbb{Q})[\hat{\phi}]/\phi(\Sha(E/\mathbb{Q})[3])$ (see also Theorem $3$ of \cite{fisherctp}). Therefore, it follows that $\mathrm{dim}_{\mathbb{F}_3}\Sha(\widehat{E}/\mathbb{Q})[\hat{\phi}]/\phi(\Sha(E/\mathbb{Q})[3])$ is even.

 Let $$2\alpha:= \mathrm{dim}_{\mathbb{F}_3}\Sha(\widehat{E}/\mathbb{Q})[\hat{\phi}]/\phi(\Sha(E/\mathbb{Q})[3]).$$ Since $\frac{|\mathrm{Sel}^{(\phi)}(E/\mathbb{Q})|}{|\mathrm{Sel}^{(\hat{\phi})}(\widehat{E}/\mathbb{Q})|}=1$ by hypothesis, we get that 
$$\beta:=\mathrm{dim}_{\mathbb{F}_3}\mathrm{Sel}^{(\phi)}(E/\mathbb{Q})=\mathrm{dim}_{\mathbb{F}_3}\mathrm{Sel}^{(\hat{\phi})}(\widehat{E}/\mathbb{Q}) .$$ If $\gamma:=\mathrm{dim}_{\mathbb{F}_3}\mathrm{Sel}^3 (E/\mathbb{Q}) $, then by the exact sequence $($\ref{exactsequencephiselmergroup}$)$ we obtain that $$\beta-\gamma+\beta-2\alpha=0$$ and, therefore, $\gamma$ must be even. There is a short exact sequence $$0 \longrightarrow{} E(\mathbb{Q})/3E(\mathbb{Q}) \longrightarrow{} \mathrm{Sel}^{3}(E/\mathbb{Q}) \longrightarrow{} \Sha(E/\mathbb{Q})[3] \longrightarrow{} 0$$ (see Theorem X.$4.2$ of \cite{aec}). Since $E/\mathbb{Q}$ has a $\mathbb{Q}$-rational point of order three and rank $0$, we obtain that $E(\mathbb{Q})/3E(\mathbb{Q}) \cong \mathbb{Z}/3\mathbb{Z}$ and, hence, $\gamma \geq 1$. Since $\gamma$ is even, we obtain that $\gamma \geq 2$. Since $\gamma=\mathrm{dim}_{\mathbb{F}_3}\mathrm{Sel}^3 (E/\mathbb{Q}) \geq 2$ and $E(\mathbb{Q})/3E(\mathbb{Q}) \cong \mathbb{Z}/3\mathbb{Z}$, $\Sha(E/\mathbb{Q})[3]$ is not trivial. Finally, since the order of $\Sha(E/\mathbb{Q})$ is a square (see \cite{aec} Corollary $17.2.1$) if $\Sha(E/\mathbb{Q})[3]$ has a nontrivial element, then $9$ divides $|\Sha(E/\mathbb{Q})|$. This concludes the proof of Part $(ii)$.

\end{proof}

\begin{emptyremark}\label{rootnumbers}
We now collect some facts concerning root numbers. Let $w(E):=\displaystyle \prod_{p \in M_{\mathbb{Q}} } w_p(E)$ where $w_p(E)$ is the local root number of $E/\mathbb{Q}$ at $p$ and $M_{\mathbb{Q}}$ is the set of places of $\mathbb{Q}$. The number $w(E)$ is called the global root number of $E/\mathbb{Q}$.
By Theorem $3.22$ of \cite{darmonmodularellipticcurves}, $\Sha(E/\mathbb{Q})$ is finite and the rank of $E/\mathbb{Q}$ is $0$. Consequently, by Theorem $1.4$ of \cite{dd1} we get that $1=(-1)^{rk(E/\mathbb{Q})}=w(E)$. 

Recall that we assume that $j_E \neq 0, 1728$. The local root number of $E/\mathbb{Q}$ at $p$ is as follows (see \cite{connel} page $95$ and \cite{rohrlich} page $132$ for $j_E \neq 0, 1728$ and $p \geq 5$, and \cite{halberstadt} page 1051 for $p=2$ or $3$)
\[ \begin{cases} 
      w_{\infty}(E)=-1 \\
      w_p(E)=1 & \mathrm{if} \ E/\mathbb{Q} \ \mathrm{has} \ \mathrm{modulo} \ p \ \mathrm{good} \ \mathrm{or} \ \mathrm{nonsplit} \ \mathrm{multiplicative} \ \mathrm{reduction}.\\
      w_p(E)=-1 & \mathrm{if} \ E/\mathbb{Q} \ \mathrm{has} \ \mathrm{modulo} \ p \ \mathrm{split} \ \mathrm{multiplicative} \  \mathrm{reduction}.\\
      w_3(E)=-1 & \mathrm{if} \ E/\mathbb{Q} \ \mathrm{has} \ \mathrm{modulo} \ 3 \ \mathrm{reduction} \ \mathrm{of} \ \mathrm{type} \ \mathrm{I}_n^*\ ,n\geq 0.
   \end{cases}
\]
\end{emptyremark}

\begin{lemma}\label{lemmatamagawanumber}
Let $E/\mathbb{Q}$ be an elliptic curve given by a Weierstrass equation of the form $($\ref{eq:3torsion}$)$ such that $E/\mathbb{Q}$ has modulo $3$ reduction of type I$_n^*$, for some $n \geq 0$. If $9 \nmid \prod_{p}c_p(E)$, then either $b=1$ or $b=r^m$ for some prime $r \neq 3$.
\end{lemma}
\begin{proof}
Since $E/\mathbb{Q}$ has modulo $3$ reduction of type I$_n^*$, for some $n \geq 0$, Proposition \ref{prop:3torsionreduction} implies that $3 \nmid b$, as we now explain. First, if $3 \mid b$ and $3 \nmid a$, then Proposition \ref{prop:3torsionreduction}, Part $(i)$, implies that $E/\mathbb{Q}$ has modulo $3$ multiplicative reduction, which is a contradiction. Moreover, if $3 \mid b$ and $3 \mid a$, then, $3^3 \nmid b$ so $3\textrm{ord}_3(a)>\textrm{ord}_3(b)$. Therefore, using Proposition \ref{prop:3torsionreduction}, Part $(ii)$, we obtain that $E/\mathbb{Q}$ has modulo $3$ reduction of type IV or IV$^*$, which is again a contradiction. 

If a prime $r$ divides $b$, then either $r^3 \nmid b$ or $r \nmid a$. Therefore, since $r \neq 3$, by Proposition \ref{prop:3torsionreduction}, Parts $(i)$ and $(ii)$, $E/\mathbb{Q}$ has reduction of type IV, IV$^*$ or (split) I$_{3\text{ord}_{r}(b)}$, and in any of these cases $3 \mid c_r(E)$. Hence, if $9 \nmid \prod_{p}c_p(E)$, then $b$ can have at most one prime divisor.
\end{proof}
\begin{proof}[\it Proof of Theorem \ref{3torsionshaandtamagawa} Part $(a)$.]
Assume that $E/\mathbb{Q}$ is semi-stable away from $3$.
If $9$ divides $\prod_{p}c_p(E)$, then Theorem \ref{3torsionshaandtamagawa}, Part $(a)$ is true. Therefore, we can assume from now on that $9 \nmid \prod_{p}c_p(E)$ and, hence, by Lemma \ref{lemmatamagawanumber} that either $b=1$ or $b=r^m$, for some prime $r$ with $r \neq 3$. We now split the proof into two cases depending on whether $b=1$ or $b=r^m$.

\medskip \underline{Case 1:} $b=1$. By a proposition of Hadano (see Theorem $1.1$ of \cite{had}) since $b=1$, $\widehat{E}/\mathbb{Q}$ is given by $$y^2+(a+6)xy+(a^2+3a+9)y=x^3.$$ The discriminant of this equation is $\widehat{\Delta}=(a^3 - 27)^3$ and $\widehat{c_4}=a(a^3+216)$.

We first show that $\widehat{E}/\mathbb{Q}$ has a prime $q$  of split multiplicative reduction. By Proposition \ref{prop:3torsionreduction}, Part $(iii)$, we obtain that $3 \mid a$ because $E/\mathbb{Q}$ has reduction of type I$_n^*$ modulo $3$. Therefore, we can write $a=3a'$ for some integer $a'$. Then $$a^2+3a+9=9((a')^2+a'+1)$$ and $(a')^2+a'+1$ cannot be a power of $3$ unless $a'=1,-2$ because it is always non zero modulo $9$. If $(a')^2+a'+1=3$, then $a'=1$ or $-2$. Moreover, $a'=1$ gives $a=3$, and together with $b=1$ gives an equation of the form $($\ref{eq:3torsion}$)$ which has discriminant $0$ and, hence, is not an elliptic curve. Therefore, $a'=1$ is not allowed. If $a'=-2$, which implies $a=-6$, we obtain, combined with $b=1$, an elliptic curve $E/\mathbb{Q}$ of the form $($\ref{eq:3torsion}$)$ which does not have reduction of type I$_n^*$ at 3 because $a^3-27=(-6)^3-27=243$ is not divisible by $3^6$, see Proposition \ref{prop:3torsionreduction}, Part $(iii)$. Moreover, if $(a')^2+a'+1= \pm 1$, then $a'=-1$.  If $a'=-1$, which implies $a=-3$, we obtain, combined with $b=1$, an elliptic curve $E/\mathbb{Q}$ of the form $($\ref{eq:3torsion}$)$ which does not have reduction of type I$_n^*$ at 3 because $a^3-27=(-3)^3-27=-54$ is not divisible by $3^6$, see Proposition \ref{prop:3torsionreduction}, Part $(iii)$. What this proves is that $a^2+3a+9$ has a prime divisor not equal to $3$, say $q$. By looking at the Weierstrass equation for $\widehat{E}/\mathbb{Q}$ we get that $\widehat{E}/\mathbb{Q}$ has split multiplicative reduction modulo $q$. Indeed, if $q \mid a^2+3a+9$ and $q \mid a+6$, then $$q \mid a^2+3a+9-(a+6)^2=a^2+3a+9-a^2-12a-36=-9(a+3)$$ and, hence, $q \mid a+3$. Since $q \mid a+6$, we obtain that $q=3$, which is a contradiction. This proves that $q \nmid a+6$ and, thus, $\widehat{E}/\mathbb{Q}$ has split multiplicative reduction modulo $q$. 

Since $\widehat{E}/\mathbb{Q}$ has split multiplicative modulo $q$, we get that $c_q(\widehat{E})=\mathrm{ord}_q((a^3-27)^3)=3\mathrm{ord}_q(a^3-27)$. Since $E/\mathbb{Q}$ and $\widehat{E}/\mathbb{Q}$ are isogenous over $\mathbb{Q}$, they have
the same $L$-function (see Korollar 1 of \cite{faltingsisogeny}) and, hence, $a_q(E)=a_q(\widehat{E})$ (see page $366$ of \cite{diamondshurman}). Moreover, $a_q(\widehat{E})=1$ because $\widehat{E}/\mathbb{Q}$ has split multiplicative reduction modulo $q$ (see page $329$ of \cite{diamondshurman}), which implies that $a_q(E)=1$. Therefore, $\widehat{E}/\mathbb{Q}$ has split multiplicative reduction modulo $q$ and $c_q(E)=\mathrm{ord}_q(\Delta)$. However, $\Delta=a^3-27$ so $c_q(E) = \mathrm{ord}_q(a^3-27)$, which implies that $\textrm{ord}_3 (\frac{c_q(\widehat{E})}{c_q(E)})>0$.

 We will now show that $9$ divides $|\Sha(E/\mathbb{Q})$|. By Part $(i)$ of Lemma \ref{tamagawagreaterorequalto2}, it is enough to show that $\text{ord}_3(\frac{\prod_{p}c_p(\widehat{E})}{\prod_{p}c_p(E)})\geq 2.$ We achieve this in the following claim.
\begin{claim}
We have $$\text{ord}_3(\frac{\prod_{p}c_p(\widehat{E})}{\prod_{p}c_p(E)}) \geq 2.$$
\end{claim}
\begin{proof}

First, we claim that because the analytic rank of $E/\mathbb{Q}$ is zero, $E/\mathbb{Q}$ has an even number of places of split multiplicative reduction. Indeed, by \ref{rootnumbers} $w(E)=1$, $w_3(E)=-1$, and $w_{\infty}(E)=-1$. Moreover, $E/\mathbb{Q}$ is semi-stable away from $3$ and by \ref{rootnumbers} for $p\neq 3$ we obtain that $w_p(E)=-1$ if and only if $E/\mathbb{Q}$ has split multiplicative reduction modulo $p$. This proves that $E/\mathbb{Q}$ has an even number of places of split multiplicative reduction.

 Let $p \neq 3$ be any prime such that $E/\mathbb{Q}$ has multiplicative reduction modulo $p$. If $E/\mathbb{Q}$ has nonsplit multiplicative reduction modulo $p$, then by lines $5,6$, and $7$ of Theorem $6.1$ of \cite{dd} we obtain that $\textrm{ord}_3(\frac{c_p(\widehat{E})}{c_p(E)})=0$. Note that in Theorem $6.1$ of \cite{dd}, $\widehat{E}$ is denoted by $E'$ and $\delta,\delta'$ are the valuations of the two discriminants. Assume now that $E/\mathbb{Q}$ has split multiplicative reduction modulo $p$. Recall that we assume that $b=1$ so $\Delta=a^3-27$ and $\widehat{\Delta}=(a^3-27)^3$. If $\textrm{ord}_p(\Delta)=\gamma$, then $\textrm{ord}_p(\widehat{\Delta})=3\gamma$. Therefore, by line $3$ of Theorem $6.1$ of \cite{dd} we obtain that $\textrm{ord}_3(\frac{c_p(\widehat{E})}{c_p(E)})=1$. Since $E/\mathbb{Q}$ is semi-stable away from $3$, the above arguments prove that if $p \neq 3$ is a prime, then $E/\mathbb{Q}$ has split multiplicative reduction modulo $p$ if and only if $\textrm{ord}_3(\frac{c_p(\widehat{E})}{c_p(E)})=1$, and $\textrm{ord}_3(\frac{c_p(\widehat{E})}{c_p(E)})=0$ otherwise.
 
 We now claim that ord$_3(\frac{\prod_{p}c_p(\widehat{E})}{\prod_{p}c_p(E)})$ is even. By line $13$ of Theorem $6.1$ of \cite{dd}, since $E/\mathbb{Q}$ has modulo 3 reduction of type I$_n^*$ for some $n\geq 0$, we obtain that $\textrm{ord}_3 (\frac{c_3(\widehat{E})}{c_3(E)})=0$.  Combining all the above with the fact that $E/\mathbb{Q}$ has an even number of places of split multiplicative reduction, we obtain that ord$_3(\frac{\prod_{p}c_p(\widehat{E})}{\prod_{p}c_p(E)})$ is even.

 Finally, since $\textrm{ord}_3 (\frac{c_q(\widehat{E})}{c_q(E)})>0$ and ord$_3(\frac{\prod_{p}c_p(\widehat{E})}{\prod_{p}c_p(E)})$ is even, we get that $\text{ord}_3(\frac{\prod_{p}c_p(\widehat{E})}{\prod_{p}c_p(E)}) \geq 2$. This proves our claim.
\end{proof}

\underline{Case 2:} $b=r^m$ for some prime $r$, and $m>0$. Recall that we assume that $E/\mathbb{Q}$ has semi-stable reduction outside of $3$ and reduction of type I$_n^*$ modulo $3$, for some $n \geq 0$. Lemma \ref{lemmatamagawanumber} implies that $r \neq 3$.

We claim that $E/\mathbb{Q}$ has split multiplicative reduction modulo $r$. Indeed, by our assumption $E/\mathbb{Q}$ is semi-stable away from $3$ and $r \neq 3$. Since $\textrm{ord}_r(b)>0$ and we know that either $r \nmid a$ or $r^3 \nmid b$, we obtain that either $3\textrm{ord}_r(a)>\textrm{ord}_r(b)$ or $\textrm{ord}_r(a)=0$. The case $3\textrm{ord}_r(a)>\textrm{ord}_r(b)$ gives that $E/\mathbb{Q}$ has reduction of type IV or IV$^*$ modulo $r$, by  Part $(ii)$ of Proposition \ref{prop:3torsionreduction}, which contradicts our assumption that $E/\mathbb{Q}$ has semi-stable reduction outside of $3$. If $\textrm{ord}_r(a)=0$, then we obtain that $3\textrm{ord}_r(a)<\textrm{ord}_r(b)$ and, hence, by Proposition \ref{prop:3torsionreduction}, Part $(i)$ we get that $E/\mathbb{Q}$ has split multiplicative reduction modulo $r$. 

If $3 \mid m$, then by Proposition \ref{prop:3torsionreduction}, Part $(i)$ we obtain that $3m \mid \prod_{p}c_p(E)$ and, hence, $9 \mid \prod_{p}c_p(E)$. Thus, if $3 \mid m$, then Part $(a)$ of Theorem \ref{3torsionshaandtamagawa} is satisfied. Therefore, we can assume that $3 \nmid m$ in what follows. Finally, since $3\textrm{ord}_r(a)<\textrm{ord}_r(b)$, by Theorem $4.1$ of \cite{koz} we obtain $\textrm{ord}_3(\frac{c_r(\widehat{E})}{c_r(E)})=-1$.

\begin{claim}\label{claimevenorder1}
The number ord$_3(\frac{\prod_{p}c_p(\widehat{E})}{\prod_{p}c_p(E)})$ is even and non-negative.
\end{claim}
\begin{proof}
There exists at least one prime $q$ different from $r$ such that $E/\mathbb{Q}$ has split multiplicative reduction modulo $q$. Indeed, by \ref{rootnumbers} $w(E)=1$, $w_r(E)=-1$, $w_3(E)=-1$, and $w_\infty(E)=-1$. Moreover, since $E/\mathbb{Q}$ is semi-stable away from $3$, by \ref{rootnumbers} for $p\neq 3$, we get that $w_p(E)=-1$ if and only if $E/\mathbb{Q}$ has split multiplicative reduction modulo $p$. This proves that $E/\mathbb{Q}$ has an even number of places of split multiplicative reduction. Therefore, since $E/\mathbb{Q}$ has split multiplicative reduction modulo $r$, we obtain that there is a prime $q \neq r$, such that $E/\mathbb{Q}$ has split multiplicative reduction modulo $q$.

Let $p \neq 3$ be any prime such that $E/\mathbb{Q}$ has multiplicative reduction modulo $p$. If $E/\mathbb{Q}$ has nonsplit multiplicative reduction modulo $p$, then by lines $5,6$, or $7$ of Theorem $6.1$ of \cite{dd} we obtain that $\textrm{ord}_3(\frac{c_p(\widehat{E})}{c_p(E)})=0$.

The elliptic curve $\widehat{E}/\mathbb{Q}$ is given by a Weierstrass equation of the form
\begin{align*}
    y^2+axy-9by=x^3-(a^3+27b)b,
\end{align*}
where $a,b$ are as in $($\ref{eq:3torsion}$)$ (see \cite{koz} equation $(3.2)$). The discriminant of this Weierstrass equation is $$\widehat{\Delta}=(a^3-27b)^3b.$$

We now prove that for any prime $p \neq 3,r,$ if $E/\mathbb{Q}$ has split multiplicative reduction modulo $p$, then $\textrm{ord}_3(\frac{c_p(\widehat{E})}{c_p(E)})=1$. Assume that $p \neq r$ is a prime such that $E/\mathbb{Q}$ has split multiplicative reduction modulo $p$. Note that since $b=r^m$ and $p\neq r$, $p\nmid b$. Recall that $\Delta=b^3(a^3-27b)$ and $\widehat{\Delta}=b(a^3-27b)^3$. If $\textrm{ord}_p(\Delta)=\gamma$, then $\textrm{ord}_p(\widehat{\Delta})=3\gamma$. Therefore, by line $3$ of Theorem $6.1$ of \cite{dd} we obtain that $\textrm{ord}_3(\frac{c_p(\widehat{E})}{c_p(E)})=1$.

Recall that there exists a prime $q$ different from $r$ such that $E/\mathbb{Q}$ has split multiplicative reduction modulo $q$. We have proved so far that $\textrm{ord}_3(\frac{c_r(\widehat{E})}{c_r(E)})=-1$, $\textrm{ord}_3(\frac{c_q(\widehat{E})}{c_q(E)})=1$, and that for any prime $p \neq 3,r$, if $E/\mathbb{Q}$ has split multiplicative reduction modulo $p$, then $\textrm{ord}_3(\frac{c_p(\widehat{E})}{c_p(E)})=1$. Moreover, we have proved above that if $E/\mathbb{Q}$ has nonsplit multiplicative reduction modulo $p$, then $\textrm{ord}_3(\frac{c_p(\widehat{E})}{c_p(E)})=0$. Therefore, since $E/\mathbb{Q}$ is semistable away from $3$ and has an even number of places of split multiplicative reduction, we obtain that ord$_3(\frac{\prod_{p}c_p(\widehat{E})}{\prod_{p}c_p(E)})$ is even and non-negative, as claimed.
\end{proof}

If ord$_3(\frac{\prod_{p}c_p(\widehat{E})}{\prod_{p}c_p(E)}) \geq 2$, then Part $(i)$ of Lemma \ref{tamagawagreaterorequalto2} implies that $9$ divides $ |\Sha(E/\mathbb{Q})|$.
Otherwise, ord$_3(\frac{\prod_{p}c_p(\widehat{E})}{\prod_{p}c_p(E)}) = 0$. Since the degree of $\phi$ is $3$, by Proposition $2$ of \cite{ks} we obtain that ord$_s(\frac{\prod_{p}c_p(\widehat{E})}{\prod_{p}c_p(E)}) = 0$ for every prime $s \neq 3$. Therefore, if $\frac{\Omega(E)}{\Omega(\widehat{E})}=3$, then Equation $($\ref{eq:selmerandtamagawa}$)$ implies that $\frac{|\mathrm{Sel}^{(\phi)}(E/\mathbb{Q})|}{|\mathrm{Sel}^{(\hat{\phi})}(\widehat{E}/\mathbb{Q})|}=1$ and, hence, by Part $(ii)$ of Lemma \ref{tamagawagreaterorequalto2} we obtain that $9$ divides $|\Sha(E/\mathbb{Q})|$. If $\frac{\Omega(E)}{\Omega(\widehat{E})}=1$, then, since $b=r^m$ with $3 \nmid m$, the following lemma shows that $9$ divides $|\Sha(E/\mathbb{Q})|$. This concludes the proof of Part $(a)$ of Theorem \ref{3torsionshaandtamagawa}.
\end{proof}
\begin{lemma}\label{lemmadescent}
 Let $E/\mathbb{Q}$ be as in Theorem \ref{3torsionshaandtamagawa}, and given by a Weierstrass equation as in $($\ref{eq:3torsion}$)$. Assume moreover that there exists a prime $r \neq 3$ such that $b=r^m$, for some integer $m$ coprime to $3$. If $\mathrm{ord}_3(\frac{\prod_{p}c_p(\widehat{E})}{\prod_{p}c_p(E)})=0$ and $\frac{\Omega(E)}{\Omega(\widehat{E})}=1$, then $9$ divides $|\Sha(E/\mathbb{Q})|$.
\end{lemma}
\begin{proof}

Since $\frac{\Omega(E)}{\Omega(\widehat{E})}=1$, by Equation $($\ref{eq:selmerandtamagawa}$)$, combined with $($\ref{kernelofisogeny}$)$, we get that $$|\mathrm{Sel}^{(\phi)}(E/\mathbb{Q})|=\frac{|\mathrm{Sel}^{(\hat{\phi})}(\widehat{E}/\mathbb{Q})|\cdot 3\cdot \prod_{p}c_p(\widehat{E})}{\prod_{p}c_p(E)}.$$
Since $\mathrm{ord}_3(\frac{\prod_{p}c_p(\widehat{E})}{\prod_{p}c_p(E)})=0$, we obtain that $\mathrm{dim}_{\mathbb{F}_3}\mathrm{Sel}^{(\phi)}(E/\mathbb{Q})=\mathrm{dim}_{\mathbb{F}_3}\mathrm{Sel}^{(\hat{\phi})}(\widehat{E}/\mathbb{Q})+1$. Therefore, if we can show that $\mathrm{dim}_{\mathbb{F}_3}\mathrm{Sel}^{(\hat{\phi})}(\widehat{E}/\mathbb{Q}) \geq 1$, then we get that $\mathrm{dim}_{\mathbb{F}_3}\mathrm{Sel}^{(\phi)}(E/\mathbb{Q}) \geq 2$ and, hence, by Part $(i)$ of Lemma \ref{tamagawagreaterorequalto2} we obtain that $9$ divides $|\Sha(E/\mathbb{Q})|$.

We now show that $\mathrm{dim}_{\mathbb{F}_3}\mathrm{Sel}^{(\hat{\phi})}(\widehat{E}/\mathbb{Q}) \geq 1$. There is a short exact sequence $$0 \longrightarrow{} E(\mathbb{Q})/\hat{\phi}(\widehat{E}(\mathbb{Q})) \longrightarrow{} \mathrm{Sel}^{(\hat{\phi})}(\widehat{E}/\mathbb{Q}) \longrightarrow{} \Sha(\widehat{E}/\mathbb{Q})[\hat{\phi}] \longrightarrow{} 0$$ (see Theorem X.$4.2$ of \cite{aec} applied to $\hat{\phi}: \widehat{E} \rightarrow{} E$). Recall that $E/\mathbb{Q}$ has rank $0$ and that $E(\mathbb{Q})$ contains a point of order $3$. The rank of $\widehat{E}/\mathbb{Q}$ is $0$ because it is isogenous to $E/\mathbb{Q}$. Moreover, since $b$ is not a cube, Theorem $1.1$ of \cite{had} implies that $\widehat{E}(\mathbb{Q})$ does not contain a point of order $3$. Therefore, $E(\mathbb{Q})/\hat{\phi}(\widehat{E}(\mathbb{Q}))$ contains a point of order $3$ which injects into $\mathrm{Sel}^{(\hat{\phi})}(\widehat{E}/\mathbb{Q})$. This proves that $\text{dim}_{\mathbb{F}_3}(\mathrm{Sel}^{(\hat{\phi})}(\widehat{E}/\mathbb{Q})) \geq 1$.
\end{proof}

\begin{proof}[\it Proof of Theorem \ref{3torsionshaandtamagawa} Part $(b)$.]
We know that $E/\mathbb{Q}$ has a Weierstrass equation as in $($\ref{eq:3torsion}$)$. By Proposition \ref{prop:3torsionreduction} if $p \neq 3$ is any prime such that $E/\mathbb{Q}$ has additive reduction, then $p \mid b$, the reduction type modulo $p$ is IV or IV$^*$, and $3 \mid c_p(E)$. Thus, if $E/\mathbb{Q}$ has more than two places of additive reduction, then $9\mid \prod_{p \mid N}c_p(E)$.
\end{proof}

\begin{proof}[\it Proof of Theorem \ref{3torsionshaandtamagawa} Part $(c)$.]
Our strategy in this proof is to show that if $9$ does not divide $\prod_{p|N}c_p(E)$,
then 9 divides $|\Sha(E/\mathbb{Q})|$.
We know that $E/\mathbb{Q}$ has a Weierstrass equation as in $($\ref{eq:3torsion}$)$. By assumption $E/\mathbb{Q}$ has exactly two places of additive reduction, say at $3$ and $r$. 

Assume that $9 \nmid \prod_{p|N}c_p(E)$ otherwise the theorem is proved. Our assumptions force $b=r$ or $b=r^2$, and $a=3rm$ where $m$ is an integer, as we now explain. Indeed, if $b=1$, then Parts $(i)$ and $(ii)$ of Proposition \ref{prop:3torsionreduction} imply that $E/\mathbb{Q}$ is semi-stable away from $3$ and that $E/\mathbb{Q}$ has only one place of additive reduction. If $b$ had two or more prime divisors, Parts $(i)$ and $(ii)$ of Proposition \ref{prop:3torsionreduction} imply that $9$ divides $\prod_{p|N}c_p(E)$. If  $r \nmid a$, then by Proposition \ref{prop:3torsionreduction}, Part $(i)$, we obtain that $E/\mathbb{Q}$ has multiplicative reduction modulo $r$, contrary to our assumption. This proves that $r \mid a$. Since $r\mid a$ and $E/\mathbb{Q}$ has additive reduction modulo $r$, Proposition \ref{prop:3torsionreduction} implies that $r \mid b$. Moreover, since we must have that $r \nmid a$ or $r^3 \nmid b$ in Equation $($\ref{eq:3torsion}$)$, we obtain that $b$ is equal to either $r$ or $r^2$. Finally, since $E/\mathbb{Q}$ has reduction I$_n^*$ modulo $3$, Proposition \ref{prop:3torsionreduction}, Part $(iii)$ implies that $3 \mid a$.

By Proposition \ref{prop:3torsionreduction}, Part $(ii)$, $E/\mathbb{Q}$ has reduction of type IV or IV$^*$ modulo $r$. By the table on page $46$ of \cite{tatealgorithm} for $r \neq 2$ and by Tableau IV of \cite{pap} for $r=2$, we obtain then that ord$_r(\Delta)=4$ or $8$. Therefore, Proposition $2$ of \cite{rohrlich} implies that $w_r(E)=\big( \frac{-3}{r} \big)$. If $\big( \frac{-3}{r} \big)=1$, then $r \equiv 1 \; (\text{mod} \; 3)$ and if $\big( \frac{-3}{r} \big)=-1$, then $r \equiv 2 \; (\text{mod} \; 3)$. We split the proof into two cases, when $w_r(E)=1$ and $r \equiv 1 \; (\text{mod} \; 3)$, and when $w_r(E)=-1$ and $r \equiv 2 \; (\text{mod} \; 3)$.

Let us first prove the following claim.
\begin{claim}\label{claimsplitmultiplicative}
If $p$ is any prime such that $E/\mathbb{Q}$ has split multiplicative reduction modulo $p$, then $\textrm{ord}_3(\frac{c_p(\widehat{E})}{c_p(E)})=1$.
\end{claim}
\begin{proof}
Note that $p \nmid 3b$ because $E/\mathbb{Q}$ has additive reduction modulo $3$ and $r$. Recall that $\widehat{E}/\mathbb{Q}$ is given by a Weierstrass equation of the form
\begin{align*}
    y^2+axy-9by=x^3-(a^3+27b)b,
\end{align*}
where $a,b$ are as in $($\ref{eq:3torsion}$)$ (see \cite{koz} equation $(3.2)$). The discriminant of this Weierstrass equation is $$\widehat{\Delta}=(a^3-27b)^3b.$$ Since $\Delta=b^3(a^3-27b)$ and $\widehat{\Delta}=b(a^3-27b)^3$, if $\textrm{ord}_p(\Delta)=\gamma$, then $\textrm{ord}_p(\widehat{\Delta})=3\gamma$. Therefore, by line $3$ of Theorem $6.1$ of \cite{dd} we obtain that $\textrm{ord}_3(\frac{c_p(\widehat{E})}{c_p(E)})=1$. This proves that if $p$ be any prime such that $E/\mathbb{Q}$ has split multiplicative reduction modulo $p$, then $\textrm{ord}_3(\frac{c_p(\widehat{E})}{c_p(E)})=1$.
\end{proof}

\underline{Case 1:} $w_r(E)=1$ and $r \equiv 1 \; (\text{mod} \; 3)$.

\begin{claim}
The number
$\text{ord}_3(\frac{\prod_{p}c_p(\widehat{E})}{\prod_{p}c_p(E)})$ is even and non-negative.
\end{claim}
\begin{proof}
First, \ref{rootnumbers} implies that $w_\infty(E)=-1$ and that $w_3(E)=-1$. Recall that $E/\mathbb{Q}$ is semi-stable away from $3$ and $r$. For $p\neq 3$ we have that $w_p(E)=-1$ if and only if $E/\mathbb{Q}$ has split multiplicative reduction modulo $p$. Since $w(E)=1$, we obtain that $E/\mathbb{Q}$ has an even number of places of split multiplicative reduction.

Since $r \equiv 1 \; (\text{mod} \; 3)$, we obtain that $\zeta_3 \in \mathbb{Q}_r$ and since we also have that $E/\mathbb{Q}$ has reduction of type IV or IV$^*$ modulo $r$, by line $10$ of Theorem $6.1$ of \cite{dd} we obtain that $\textrm{ord}_3(\frac{c_r(\widehat{E})}{c_r(E)})=0$. If $E/\mathbb{Q}$ has nonsplit multiplicative reduction modulo $p$, then $\textrm{ord}_3(\frac{c_p(\widehat{E})}{c_p(E)})=0$ by line $4$ of Theorem $6.1$ of \cite{dd}. Also, If $E/\mathbb{Q}$ has split multiplicative reduction modulo $p$, then $\textrm{ord}_3(\frac{c_p(\widehat{E})}{c_p(E)})=1$ by Claim \ref{claimsplitmultiplicative}. Finally, line $13$ of Theorem $6.1$ of \cite{dd} implies that $\textrm{ord}_3(\frac{c_3(\widehat{E})}{c_3(E)})=0$. Putting all those together we obtain that $\text{ord}_3(\frac{\prod_{p}c_p(\widehat{E})}{\prod_{p}c_p(E)})$ is even and greater than or equal to zero. This concludes the proof of the claim.
\end{proof}

We are now ready to conclude the proof of Case $1$. Assume first that $\text{ord}_3(\frac{\prod_{p}c_p(\widehat{E})}{\prod_{p}c_p(E)})=0.$ Since the degree of $\phi$ is $3$, we can use Proposition $2$ of \cite{ks} and obtain that ord$_p(\frac{\prod_{p}c_p(\widehat{E})}{\prod_{p}c_p(E)}) = 0$ for every prime $p \neq 3$. Hence, if $\frac{\Omega(E)}{\Omega(\widehat{E})}=3$, then Equation $($\ref{eq:selmerandtamagawa}$)$ implies that $\frac{|\mathrm{Sel}^{(\phi)}(E/\mathbb{Q})|}{|\mathrm{Sel}^{(\hat{\phi})}(\widehat{E}/\mathbb{Q})|}=1.$ Therefore, Part $(ii)$ of Lemma \ref{tamagawagreaterorequalto2} implies that $9$ divides  $|\Sha(E/\mathbb{Q})$|. If $\frac{\Omega(E)}{\Omega(\widehat{E})}=1$, then, since $b$ is equal to either $r$ or $r^2$, and we assume that $\text{ord}_3(\frac{\prod_{p}c_p(\widehat{E})}{\prod_{p}c_p(E)})=0$, Lemma \ref{lemmadescent} implies that $9$ divides $|\Sha(E/\mathbb{Q})|$. Assume now that $\text{ord}_3(\frac{\prod_{p}c_p(\widehat{E})}{\prod_{p}c_p(E)}) \geq 2.$
Then $9$ divides $|\Sha(E/\mathbb{Q})|$ by Part $(i)$ of Lemma \ref{tamagawagreaterorequalto2}. This proves  Theorem \ref{3torsionshaandtamagawa}, Part $(c)$, in the case where $w_r(E)=1$.

\medskip \underline{Case 2:} $w_r(E)=-1$ and $r \equiv 2 \; (\text{mod} \; 3)$.

\begin{claim}
The number ord$_3(\frac{\prod_{p}c_p(\widehat{E})}{\prod_{p}c_p(E)})$ is even and non-negative.
\end{claim}
\noindent {\it Proof.}
By \ref{rootnumbers} $w(E)=1$, $w_3(E)=-1$, and $w_{\infty}(E)=-1$. If $p \neq 3,r$ is a prime such that $E/\mathbb{Q}$ has split multiplicative reduction modulo $p$, then $w_p(E)=-1$.  If $p \neq 3,r$ is a prime such that $E/\mathbb{Q}$ has nonsplit multiplicative or good reduction modulo $p$, then $w_p(E)=1$. Therefore, since $E/\mathbb{Q}$ is semi-stable away from $3$ and $r$, we obtain that $E/\mathbb{Q}$ has an odd number of primes of split multiplicative reduction. Note that in particular $E/\mathbb{Q}$ has at least one prime of split multiplicative reduction.

If $p$ is any prime such that $E/\mathbb{Q}$ has split multiplicative reduction modulo $p$, then Claim \ref{claimsplitmultiplicative} implies that $\textrm{ord}_3(\frac{c_p(\widehat{E})}{c_p(E)})=1$. By line $4$ of Theorem $6.1$ of \cite{dd}, if $E/\mathbb{Q}$ has nonsplit multiplicative reduction modulo $p$, then $\textrm{ord}_3(\frac{c_p(\widehat{E})}{c_p(E)})=0$. Moreover, by line $10$ of Theorem $6.1$ of \cite{dd} we obtain that $\textrm{ord}_3(\frac{c_r(\widehat{E})}{c_r(E)})=-1$ because $r \equiv 2 \; (\text{mod} \; 3)$ implies $\zeta_3 \notin \mathbb{Q}_r$ and $E/\mathbb{Q}$ has reduction of type IV or IV$^*$ modulo $r$. Finally, line $13$ Theorem $6.1$ of \cite{dd} implies that $\textrm{ord}_3(\frac{c_3(\widehat{E})}{c_3(E)})=0$. Since there is an odd number of primes of split multiplicative reduction, and in particular at least one prime of split multiplicative reduction, we obtain that ord$_3(\frac{\prod_{p}c_p(\widehat{E})}{\prod_{p}c_p(E)})$ is even and non-negative. This proves our claim.
\end{proof}

We are now ready to conclude the proof of Case $2$. Assume first that $\text{ord}_3(\frac{\prod_{p}c_p(\widehat{E})}{\prod_{p}c_p(E)})=0.$ Since the degree of $\phi$ is $3$, by line $5$ Theorem $6.1$ of \cite{dd} (or Proposition $2$ of \cite{ks}) we obtain that ord$_q(\frac{\prod_{p}c_p(\widehat{E})}{\prod_{p}c_p(E)}) = 0$ for every prime $q \neq 3$. Therefore, if $\frac{\Omega(E)}{\Omega(\widehat{E})}=3$, Equation $($\ref{eq:selmerandtamagawa}$)$ implies that $\frac{|\mathrm{Sel}^{(\phi)}(E/\mathbb{Q})|}{|\mathrm{Sel}^{(\hat{\phi})}(\widehat{E}/\mathbb{Q})|}=1$ and, hence, by Part $(ii)$ of Lemma \ref{tamagawagreaterorequalto2} we obtain that $9$ divides $|\Sha(E/\mathbb{Q})|$. If $\frac{\Omega(E)}{\Omega(\widehat{E})}=1$, since $b$ is equal to $r$ or $r^2$, and we assume that $\text{ord}_3(\frac{\prod_{p}c_p(\widehat{E})}{\prod_{p}c_p(E)})=0$, Lemma \ref{lemmadescent} implies that $9$ divides $|\Sha(E/\mathbb{Q})|$. Assume now that $\text{ord}_3(\frac{\prod_{p}c_p(\widehat{E})}{\prod_{p}c_p(E)}) \geq 2.$
Then $9$ divides $|\Sha(E/\mathbb{Q})|$ by Part $(i)$ of Lemma \ref{tamagawagreaterorequalto2}. This proves Theorem \ref{3torsionshaandtamagawa}, Part $(c)$, in the case $w_r(E)=-1$. This completes the proof of Theorem \ref{3torsionshaandtamagawa}.
\end{proof}

\begin{remark}
The condition that $E/\mathbb{Q}$ has modulo $3$ reduction I$_n^*$ for some $n\geq0$ in Theorem \ref{3torsionshaandtamagawa} is necessary. Indeed, the conclusion of Theorem \ref{3torsionshaandtamagawa} for the elliptic curve $E/\mathbb{Q}$ with Cremona \cite{cremonabook} label 324d1 (or LMFDB \cite{lmfdb} label 324.b1) is not true. This curve has reduction II modulo $3$, reduction IV$^*$ modulo $2$, $\prod_{p}c_p(E)=3$, $\Sha(E/\mathbb{Q})$ is trivial, and $E(\mathbb{Q})\cong \mathbb{Z}/3\mathbb{Z}$.

Moreover, the condition that $E/\mathbb{Q}$ has analytic rank $0$ in Theorem \ref{3torsionshaandtamagawa} is necessary. Indeed, the conclusion of Theorem \ref{3torsionshaandtamagawa} for the elliptic curve $E/\mathbb{Q}$ with Cremona \cite{cremonabook} label 171b2 (or LMFDB \cite{lmfdb} label 171.b2) is not true. This curve has reduction I$_0^*$ modulo $3$, rank $1$, $\prod_{p}c_p(E)=6$, $\Sha(E/\mathbb{Q})$ is trivial, and $E(\mathbb{Q}) \cong \mathbb{Z} \times \mathbb{Z}/3\mathbb{Z}$.
\end{remark}

\begin{remark}\label{remarkagashestein}
Let $E/\mathbb{Q}$ be an optimal elliptic curve of analytic rank $0$. It follows from work of Agashe and Stein that the Birch and Swinnerton-Dyer conjecture, combined with the conjecture that the Manin constant is $1$, imply that {\it the odd part of $|E(\mathbb{Q})|$ divides $|\Sha(E/\mathbb{Q})| \cdot \prod_p c_p(E)$} (see the end of Section 4.3 of \cite{agashestein}). Lorenzini has proved the above statement up to a power of $3$ (see Proposition 4.2 of \cite{lor}). 

Without the assumption that $E/\mathbb{Q}$ is optimal, if $E/\mathbb{Q}$ has reduction of type I$_n^*$ modulo $3$, then Theorem \ref{3torsionshaandtamagawa} proves that the odd part of $|E(\mathbb{Q})|$ divides $|\Sha(E/\mathbb{Q})| \cdot \prod_p c_p(E)$. Note however that the curve $E/\mathbb{Q}$ with Cremona \cite{cremonabook} label 14a4 has a $\mathbb{Q}$-rational point of order $3$, $\prod_p c_p(E)=2$, and $|\Sha(E/\mathbb{Q})|=1$. Thus the assumption that the curve $E/\mathbb{Q}$ is optimal is necessary for the statement that the odd part of $|E(\mathbb{Q})|$ divides $|\Sha(E/\mathbb{Q})| \cdot \prod_p c_p(E)$ to be true in general.
\end{remark}
 
\begin{theorem}\label{prop:dpowerof3}
Let $E/\mathbb{Q}$ be an elliptic curve of conductor $N$ with $3 \nmid N$. Let $d = - 3$ or $3$ and assume that $L(E^d,1) \neq 0$. Then
\begin{enumerate}
  \item If $d=-3$, then $
 |E^d(\mathbb{Q})|^2 \quad \text{divides} \quad |\Sha(E^d/\mathbb{Q})| \cdot \displaystyle\prod_{p|N}c_p(E^d), \quad \text{up to a power of } 2.$
     \item If $d=3$, then $
 |E^d(\mathbb{Q})|^2 \quad \text{divides} \quad |\Sha(E^d/\mathbb{Q})| \cdot \displaystyle\prod_{p|2N}c_p(E^d), \quad \text{up to a power of } 2.$
 \end{enumerate}
\end{theorem}

\begin{proof}
 We first claim that
 $|E^d(\mathbb{Q})|^2$ divides $|\Sha(E^d/\mathbb{Q})| \cdot \prod_{p}c_p(E^d)$, up to a power of 2. By Corollary \ref{corollarytorsion}, the only primes that can divide $|E^d(\mathbb{Q})|$ are $2$ and $3$. Therefore, if $E^d/\mathbb{Q}$ does not have a point  of order $3$, our claim is proved since $E^d(\mathbb{Q})$ has order a power of $2$.  So assume that $E^d(\mathbb{Q})$ contains a $\mathbb{Q}$-rational point of order $3$. A theorem of Lorenzini (see \cite{lor} Proposition $1.1$) implies that if $E/\mathbb{Q}$ is an elliptic curve with a $\mathbb{Q}$-rational point of order $9$, then $9^2$ divides $\prod_{p} c_p(E)$ except for the curve that has Cremona label \cite{cremonabook} 54b3 with $\prod_{p} c_p(E)=27$. If $E/\mathbb{Q}$ is an elliptic curve with $3 \nmid N$, then $E/\mathbb{Q}$ has good reduction I$_0$ modulo $3$. Therefore, if $d=-3$ or $3$, we obtain, using Proposition $1$ of \cite{com}, that $E^d/\mathbb{Q}$ has reduction of type I$_0^*$ modulo $3$. However,  the curve with Cremona label \cite{cremonabook} 54b3 does not have reduction of type I$_0^*$ modulo $3$ and, hence, it cannot be a twist by $d$ of an elliptic curve with good reduction modulo $3$. Therefore, we can assume that $E^d(\mathbb{Q})$ does not contain a point of order $9$.
 
 Since we assume that $E^d(\mathbb{Q})$ does not contain a point of order $9$ and the only primes that can divide $|E^d(\mathbb{Q})|$ are $2$ and $3$, to prove that
 $|E^d(\mathbb{Q})|^2$ divides $|\Sha(E^d/\mathbb{Q})| \cdot \prod_{p}c_p(E^d)$, up to a power of 2, it is enough to show that $9$ divides $|\Sha(E^d/\mathbb{Q})| \cdot \prod_{p}c_p(E^d)$. Since $3 \nmid N$ and, hence, $E/\mathbb{Q}$ has good reduction I$_0$, by Proposition $1$ of \cite{com}, we obtain that $E^d/\mathbb{Q}$ has reduction of type I$_0^*$ modulo $3$. If $E^d/\mathbb{Q}$ is semi-stable away from $3$, then applying Part $(a)$ of Theorem \ref{3torsionshaandtamagawa} to $E^d/\mathbb{Q}$ proves that $9$ divides $|\Sha(E^d/\mathbb{Q})| \cdot \prod_{p}c_p(E^d)$. If $E^d/\mathbb{Q}$ has more than two places of additive reduction, then by Part $(b)$ of Theorem \ref{3torsionshaandtamagawa} we obtain that $9 \mid \prod_{p}c_p(E^d)$. Finally, if $E^d/\mathbb{Q}$ has exactly two places of additive reduction, then Part $(c)$ of Theorem \ref{3torsionshaandtamagawa} implies that $9$ divides $|\Sha(E^d/\mathbb{Q})| \cdot \prod_{p}c_p(E^d)$.

\textit{Proof of $(i)$:} We just need to show that the product can be taken over all the primes of bad reduction for $E/\mathbb{Q}$. If $d=-3$, then the only prime that ramifies in $\mathbb{Q}(\sqrt{d})$ is $3$. Therefore the primes of bad reduction of $E^d/\mathbb{Q}$ consist of the primes of bad reduction of $E/\mathbb{Q}$ as well as 3. Hence, it is enough to show that $3 \nmid c_3(E^d)$. By Proposition $1$ of \cite{com} we see that $E^d/\mathbb{Q}$ has modulo $3$ reduction of type I$_0^*$ and $c_3(E^d)=1, 2$ or $4$ by page $367$ of \cite{silverman2}.

\textit{Proof of $(ii)$:} If $d=3$, then $2$ and $3$ are the only primes that ramify in $\mathbb{Q}(\sqrt{d})$. Therefore the primes of bad reduction of $E^d/\mathbb{Q}$ consist of all the odd primes of bad reduction of $E/\mathbb{Q}$ as well as $3$ and possibly $2$. By Proposition $1$ of \cite{com} we see that $E^d/\mathbb{Q}$ has modulo $3$ reduction of type I$_0^*$ and $c_3(E^d)=1, 2$ or $4$ by page $367$ of \cite{silverman2}. This completes the proof.
\end{proof}

We are now ready to prove a slightly stronger form of Agashe's Conjecture \ref{agasheconjecture}.

\begin{corollary}\label{corollaryagasheconjecture}
Let $E/\mathbb{Q}$ be an elliptic curve of conductor $N$ and let $-D$ be a negative fundamental discriminant such that $D$ is coprime to $N$. Suppose that $L(E^{-D},1) \neq 0$. Then
\vspace{-.1cm}
\begin{align*}
 |E^{-D}(\mathbb{Q})|^2 \quad \text{divides} \quad |\Sha(E^{-D}/\mathbb{Q})| \cdot \prod_{p|N}c_p(E^{-D}), \quad \text{up to a power of 2.} 
\end{align*}
\end{corollary}
\begin{proof}
Recall that an integer is called a fundamental discriminant if it is the discriminant of a quadratic number field. If $n$ is a square-free integer, then the discriminant of $\mathbb{Q}(\sqrt{n})$ is $n$ in the case where $n \equiv 1 \; (\text{mod}\; 4)$ and $4n$ otherwise. If $-D$ is square-free, which happens in the case where $-D \equiv 1 \; (\text{mod} \; 4)$, then Part $(ii)$ of Corollary \ref{corollarytorsion} implies that the only prime that can divide $|E^{-D}(\mathbb{Q})_{tors}|$ is $2$, except possibly for the case where $D=3$. We note that since $-1$ is not a fundamental discriminant, we must have that $D \neq 1$. If $D=3$, then Theorem \ref{prop:dpowerof3} implies the desired result.

Assume now that $-D$ is not square-free. Then we know that $-D=4n$, where $n$ is a square-free integer. Since $E^{-D}=E^n$ and $n$ is square-free, by applying Part $(i)$ of Corollary \ref{corollarytorsion} we obtain that the only prime that can divide $|E^{-D}(\mathbb{Q})_{tors}|$ is $2$, except possibly for the cases where $n=-1$ or $-3$.  If $n=-3$, then since $E^{-12}=E^{-3}$, Part $(i)$ of Theorem \ref{prop:dpowerof3} implies the desired result. Therefore, in order to prove our corollary it remains to provide a proof for the case where $n=-1$, i.e., where $-D=-4$. 

Assume that $D=4$. By our assumption $2 \nmid N$. We will show that $E^{-4}(\mathbb{Q})=E^{-1}(\mathbb{Q})$ can only contain points of order a power of $2$. Recall that since $L(E^{-1},1) \neq 0$, Theorem $3.22$ of \cite{darmonmodularellipticcurves} implies that the rank of $E^{-1}/\mathbb{Q}$ is zero. Since $2 \nmid N$, by Tables I and II of \cite{com} we obtain that $E^{-1}/\mathbb{Q}$ has reduction of type I$_4^*$, I$_8^*$, II, or II$^*$ modulo $2$. If  $E^{-1}/\mathbb{Q}$ had a point of order $5$ or $7$, then $E^{-1}/\mathbb{Q}$ would have had multiplicative reduction modulo $2$ which is impossible. Therefore, the only primes that can divide $|E^{-1}(\mathbb{Q})|$ are $2$ and $3$. Assume that $E^{-1}/\mathbb{Q}$ has a point of order $3$, and we will arrive at a contradiction. Since $E^{-1}/\mathbb{Q}$ has a point of order $3$, $E^{-1}/\mathbb{Q}$ has an equation of the form $($\ref{eq:3torsion}$)$. By Proposition \ref{prop:3torsionreduction}, $E^{-1}/\mathbb{Q}$ has either semi-stable reduction or reduction of type IV or IV$^*$ modulo $2$. However, this is a contradiction because $E^{-1}/\mathbb{Q}$ has reduction of type I$_4^*$, I$_8^*$, II or II$^*$ modulo $2$. This proves that $E^{-1}(\mathbb{Q})$ only contains points of order a power of $2$. This concludes our proof.
\end{proof}

\bibliographystyle{plain}
\bibliography{bibliography.bib}

\end{document}